\documentclass[11pt]{article}

\def\abs#1{\ensuremath{\left|#1\right|}}
\newtheorem{theorem}{Theorem}

\def\square{\hbox{\vrule\vbox{\hrule\phantom{o}\hrule}\vrule}}
\def\endofproofsymbol{\square}

\newenvironment{proof}{\par\noindent\textbf{Proof.}\hskip\labelsep}
                      {\unskip\ \nobreak\hbox{}\hfill
                       \endofproofsymbol{\parfillskip=0pt\par}
                       \vspace{\baselineskip}}

\def\eqref#1{\mbox{(\ref{eq:#1})}}    

\def\BC#1#2{\ensuremath{\left({#1\atop#2}\right)}}  

\def\Lfloor{\ensuremath\left\lfloor}
\def\Rfloor{\ensuremath\right\rfloor}

\def\Lceil{\ensuremath\left\lceil}
\def\Rceil{\ensuremath\right\rceil}

{\makeatletter
\@addtoreset{equation}{section}}

\def\NR{\mathit{NR}}

\begin{document}

\title{A characterization of the Frobenius problem \\ and its \\ 
       application to arithmetic progressions\footnote{An early
         version of the material in Sections~\ref{Sect:2} and~\ref{Sect:3} was
         published in working paper 11-99 of the Schulich School of
         Business in July 1999.}
       }
\author{Hans J. H. Tuenter \\
       {\footnotesize  Schulich School of Business,} \\ 
       {\footnotesize  York University, Toronto,} \\
       {\footnotesize  Canada, M3J 1P3}, \\
       {\footnotesize Email: htuenter@schulich.yorku.ca}}  
\date{\small \begin{tabular}{ll}
              First version: & July 1999  \\
              This version:  & June 2006 
            \end{tabular}      
}

\maketitle
\begin{abstract}
In the Frobenius problem we are given a set of coprime, positive integers $a_1, a_2,\ldots,a_k$,
and are interested in the set of positive numbers $\NR$ 
that have no representation by the linear form $\sum_i a_ix_i$ in nonnegative integers $x_1, x_2,\ldots,x_k$.
We give a functional relationship that completely characterizes the set $\NR$, 
and apply it to the case when the numbers are in an arithmetic progression.
%
%
\end{abstract}

\section{Introduction}
Given a set $\mathcal{A}$ of coprime, positive integers $a_1, a_2,\ldots,a_k$, 
one can define the set of positive numbers~$\NR$ that have no representation by the linear form $\sum a_ix_i$ in nonnegative integers $x_1, x_2,\ldots,x_k$.
It is well-known that $\NR$ is finite, and that there is an integer $L$, 
such that every integer larger than $L$ is {\em representable} 
as a nonnegative linear combination of $a_1, a_2,\ldots,a_k$. 

This is a classic result, attributed to Issai Schur (1875--1941) by Alfred Theodore Brauer (1894--1985),
who was his doctoral student and teaching assistant at the University of Berlin~\cite[p.~133]{Curtis:1999}.
Brauer is also the one who 
coined the phrase ``The Frobenius Problem''
as the problem of determining the largest element of $\NR$.
He named it after Georg Ferdinand Frobenius (1849--1917) who would mention it occasionally in his lectures~\cite{Brauer:1942}. Frobenius, also at the University of Berlin, was Schur's thesis advisor.

Nowadays, the scope of the Frobenius problem is usually extended to include questions relating to the entire set of integers that are not representable.
Several questions of interest, such as the cardinality of the set $\NR$, are related to the so-called Sylvester sums
\[ S_m=\sum_{n\in\NR}n^m,
\] 
where $m$ is a nonnegative integer. These sums are named 
in honour of James Joseph Sylvester (1814--1897), 
who posed the question of ``How many integers are not representable'' for the two-variable case~\cite{Sylvester:1883}.

Although the history of the Frobenius problem goes back to at least the late nineteenth century, it still is an active area
of research, as witnessed by the recent monograph by Ram{\'\i}rez Alfons{\'\i}n~\cite{RamirezAlfonsin:2006}, 
and its bibliography of several hunderd contemporary references.

In the current paper, we extend the characterization for the two-variable case, derived in~\cite{Tuenter:2006}, to the general case, apply it to the Frobenius problem for arithmetic progressions,
and derive closed-form expressions for their Sylvester sums. 

\section{A Characterization}  \label{Sect:2}
For a particular element $a$ of the set $\mathcal{A}$, 
let $\mathcal{N}$ be the set of positive integers such that, if $n$ is representable, then $n-a$ is not.
As no two elements of $\mathcal{N}$ can differ by a multiple of $a$, 
and the numbers $-a+1,\ldots,-1$ are not representable,
it follows that the set $\mathcal{N}$ is a complete, positive residue set modulo~$a$,
and hence has cardinality $\abs{\mathcal{N}}=a-1$.
The set~$\mathcal{N}$ is a known entity in the literature on the Frobenius problem,
and first appeared in a paper by Brauer and Shockley~\cite{BrauerShockley:1962}.
Note that $\NR$ and $\mathcal{N}$ are equivalent characterizations, in the sense that one uniquely determines the other.
So, it seems only natural that any characteristic of $\NR$ can be described in terms of the set $\mathcal{N}$.
This is indeed the case as evidenced by the following theorem.
%
\begin{theorem}  \label{Theo:1}
For every function $f$, the following identity holds:
\begin{equation} 
\sum_{n\in\NR}\left[\vphantom{\sum}
                       f(n+a)-f(n)\right] =
   \sum_{n\in{\cal N}} f(n) - \sum_{n=1}^{a-1} f(n).
 \label{eq:2.1}
\end{equation}
\end{theorem}
\begin{proof}
Every element of $\mathcal{N}$, can be represented as $n_j=j+am_j$, $j=1,\ldots,a-1$, with $m_j$ a nonnegative integer.
This gives
\begin{eqnarray*}
   \sum_{n\in\NR}\left[f(n+a)-f(n)\right]
  &=&
   \sum_{j=1}^{a-1}\sum_{m=0}^{m_j-1}\left[f(j+ma+a)-f(j+ma)\right], \\[0.25cm]
\noalign{\noindent 
where the inner sum is void when $m_j=0$. The inner sum is a telescoping
sum, and simplifies to}
  &=&
   \sum_{j=1}^{a-1}\left[f(j+m_ja)-f(j)\right]=\sum_{n\in{\cal N}}f(n)-\sum_{n=1}^{a-1}f(n),
\end{eqnarray*}
completing the proof.
\end{proof}
Note that the cardinality of $\mathcal{N}$ is easily derived by taking
$f(n)$ as a nonzero constant.

\section{Applications}  \label{Sect:3}
With the characterization  
at our disposal, we are now in a position to derive various properties 
of~$\NR$ in terms of the set~$\mathcal{N}$.
As a first application, take $f(n)$ as a threshold function: $f(n)=1$, for $n\ge\tau$,
and 0, otherwise, setting the threshold $\tau$, at the largest element of $\mathcal{N}$.
The right-hand side of~\eqref{2.1} evaluates to one, so that there must be an $n\in\NR$ such that $n+a\ge\tau$.
It is easily seen that increasing the threshold renders the right-hand side of~\eqref{2.1} as zero, 
so that there cannot be an element $n\in\NR$, such that $n+a>\tau$. 
This gives Lemma~3 of Brauer and Shockley~\cite{BrauerShockley:1962}:
\begin{equation}
  L+a=\max\ \{n\mid n\in\mathcal{N}\}.
  \label{eq:3.1}
\end{equation}
Theorem~2.3 in Selmer~\cite{Selmer:1977}, an expression for
the number of integers that are not representable,
is obtained by taking $f(n)=n$:
\begin{equation}
 S_0=\frac1{a}\sum_{n\in{\mathcal N}}n-\frac12(a-1)=\sum_{n\in{\mathcal N}}\Lfloor n/a\Rfloor,
 \label{eq:3.2}
\end{equation}
where the last equality follows from the above characterization of the elements of~${\mathcal N}$ as a residue set.
The latter equality can also be obtained by taking $f(n)=\Lfloor n/a\Rfloor$, 
where $\Lfloor x\Rfloor$ denotes the greatest integer that is not larger than~$x$.

As another application, take $f(n)=e^{nz}$ and extend the summation on the right-hand side of~\eqref{2.1} to include
$n=0$. Now use the sum formula for a finite geometric series, and divide by $e^{az}-1$ to give the exponential
generating formula (egf) of the set $\NR$ as
\begin{equation}
  \sum_{n\in\NR} e^{nz} = 
    \frac{1}{e^{az}-1}\sum_{n\in\mathcal{N}_0}e^{nz}-\frac{1}{e^z-1},  
  \label{eq:3.3}
\end{equation}
where $\mathcal{N}_0$ is the set $\mathcal{N}$ with the zero element added.
Note that, by letting $z$ go to infinity, we recover~\eqref{3.1}.\
By expanding the exponential in the left-hand side of~\eqref{3.3},
and changing the order of summation we see that we have also obtained 
the exponential generating function of 
the Sylvester sums: 
\begin{equation}
   \sum_{m=0}^{\infty}S_m\frac{z^m}{m!} = \sum_{n\in\NR} e^{nz}
   = \sum_{n\in\mathcal{N}_0}\frac{e^{nz}}{e^{az}-1}-\frac{1}{e^z-1}.
  \label{eq:3.1.1}
\end{equation}
In this we can recognize the footprint of the Bernoulli numbers and polynomials through their generating function:
\begin{equation}
   \frac{te^{tx}}{e^t-1}=\sum_{m=0}^\infty B_m(x)\frac{t^m}{m!},
 \label{eq:4.4.xx}
\end{equation}
where the Bernoulli numbers, $B_m$ are simply the Bernoulli polynomials evaluated at
zero: $B_m=B_m(0)$.
Multiplying both sides of~\eqref{3.1.1} by~$z$, and using~\eqref{4.4.xx} with a change of variables
gives an expression for the Sylvester sums in terms of the Bernoulli polynomials:
\begin{equation}
   mS_{m-1} 
   = a^{m-1}\sum_{n\in{\cal N}_0}B_m(n/a)-B_m.
  \label{eq:3.1.3} 
\end{equation}
This result can also be obtained directly from Theorem~\ref{Theo:1} by taking $f(n)=B_m(n/a)$, 
using the difference formula for the Bernoulli polynomials, 
and Raabe's multiplication theorem.

\section{The Frobenius problem for two variables}
In the classic Frobenius problem, as posed by Sylvester~\cite{Sylvester:1883}, the set $\mathcal{A}$ consists of
two coprime integers~$a$ and~$b$. This problem has been extensively studied,
and a summary of the results can be found in~\cite{Tuenter:2006}.
Selmer~\cite[Example~3.1]{Selmer:1977} gives the set~$\mathcal{N}$ in the
two-variable case:
\begin{equation}
  {\cal N} = \left\{\vphantom\sum bn,\ n=1,\ldots,a-1\right\},
  \label{eq:5.1.0}
\end{equation}
and applies~\eqref{3.1} and~\eqref{3.2} to give the well-known and classic results
$L=ab-a-b$ and $S_0=\frac12(a-1)(b-1)$.
Characterization~\eqref{3.1.3} gives
\begin{equation}
   mS_{m-1} 
   = a^{m-1}\sum_{n=0}^{a-1}B_m(nb/a)-B_m.
  \label{eq:4.1.x}  
\end{equation}
While this is an elegant and succinct representation, and determines the Sylvester sums as an explicit
function of the parameters $a$ and $b$, 
it is not the most convenient representation from a computational point of view
for large values of~$a$.
To arrive at such a representation, one can expand the Bernoulli polynomial, change the order of summation,
and use the well-known expression for the sum-of-powers to derive an explicit expression for~$S_m$.
This approach was taken in~\cite{Tuenter:2006}, and further details can be found there. 
Here, we take a different route and take the exponential generating function as a starting point.
The egf is easily determined using~\eqref{3.3} as
\begin{equation}
  \sum_{m=0}^\infty S_m\frac{z^m}{m!} =
  \frac{e^{abz}-1}{(e^{az}-1)(e^{bz}-1)}-\frac1{e^z-1}.
  \label{eq:5.x.1}
\end{equation}
Now multiply both sides of this equation by $z^2$, and use the
convolution property of the 
egf,\footnote{Convolution property: If the sequence $S_m$ has the egf $f(z)$, and the sequence $T_m$ the egf $g(z)$,
then $f(z)g(z)$ is the egf of the sequence $S_m\star T_m=\sum_{j=0}^m\BC{m}{j}S_{m-j}T_j$.}
with the egf of the Bernoulli
polynomials to give:
\begin{equation}
  m(m-1)S_{m-2} =
  \sum_{j=0}^m\BC{m}{j}a^{m-1-j}B_{m-j}b^{j-1}\phi_j(a)-mB_{m-1},
 \label{eq:5.x.1b}
\end{equation}
where we have introduced the Bernoullian polynomial $\phi_j(x)=B_j(x)-B_j$,
for notational brevity.\footnote{The polynomials $\phi_j(x)$ have the egf $t\frac{e^{tx}-1}{e^x-1}$, 
and this is actually an older definition of the Bernoulli
polynomials, that is no longer in common use. 
Whittaker~\cite[p.~98]{Whittaker:1902} refers to $\phi_j(x)$ as the Bernoullian polynomials.
Of course, apart from a slight notational advantage in the current paper, the difference is immaterial.
}
This gives an $O(m^2)$ computational scheme for the Sylvester sums, as opposed to an $O(am)$ scheme using~\eqref{4.1.x}.
To make the dependence upon the parameters~$a$ and~$b$ more explicit, we can expand the Bernoullian polynomials.
Since $\phi_0(x)=0$, we can drop the index
$j=0$ from the summation, replace $m$ by $m+1$, and simplify to give
the expression derived by R\"odseth~\cite{Rodseth:1994}:
\[ S_{m-1} = \frac1{m(m+1)}\sum_{i=0}^{m}\sum_{j=0}^{m-i}\BC{m+1}{i}\BC{m+1-i}{j}B_iB_ja^{m-j}b^{m-i}-\frac1mB_{m}.
\]
This gives $S_0=\frac12(a-1)(b-1)$,
$S_1=\frac1{12}(a-1)(b-1)(2ab-a-b-1)$, 
and $S_2=\frac1{12}ab(a-1)(b-1)(ab-a-b)$.
The reason to repeat these expressions here, in particular the
egf~\eqref{5.x.1}, is to showcase the simularity with the formulae
for the arithmetic case (and its generalization) that we derive in the next sections.

\section{The Frobenius problem for an arithmetic progression} \label{Sect:1}
Given the arithmetic progression $a, a+d, a+2d, \ldots, a+sd$,
where $a$ and $d$ are two positive integers that are relative prime, 
one can define the set of positive numbers~$\NR$ that have no representation by the linear form $\sum_{i=0}^s(a+id)x_i$ in nonnegative integers $x_0, x_1,\ldots,x_s$. 
The choice of $a=1$ leads to the trivial case that the set $\NR$ is empty, so we will assume that $a\ge2$. 
Furthermore, one can make the assumption $s<a$, as otherwise the last element of the arithmetic progression can be expressed as $a+sd=d\times a+1\times(a+(s-a)d)$,
and thus the set $\NR$ does not alter if we restrict ourselves to the arithmetic progression
$a, a+d,\ldots,a+(s-1)d$.
Note that, when $s<a$, there are no redundancies in 
$\mathcal{A}=\left\{a, a+d, \ldots, a+sd\right\}$, 
and none of its elements can be expressed as a nonnegative linear combination
of the others.

Roberts~\cite{Roberts:1956} showed that the largest element of $\NR$ is given by
\begin{equation}
  L+1= 
  \left(\Lfloor{a-2\over s}\Rfloor+1\right)a+(a-1)(d-1),
  \label{eq:1.2}
\end{equation}
generalizing the well-known and classic result 
for the case $s=1$,
and a Theorem by Brauer~\cite{Brauer:1942} for the case $d=1$.
A simplified proof of~\eqref{1.2} was given by Bateman~\cite{Bateman:1958}.

Grant~\cite{Grant:1973} determined the cardinality of the set $\NR$ as
\begin{equation}
   S_0 = 
   \frac12\left(\Lfloor{a-2\over s}\Rfloor+1\right)(a+t)+\frac12(a-1)(d-1),
  \label{eq:1.3}
\end{equation}
where $t= (a-2) \, \bmod s$. 
This generalizes the result by Sylvester~\cite{Sylvester:1883} for the case $s=1$,
and the result by Nijenhuis and Wilf~\cite{NijenhuisWilf:1972} for $d=1$.
The proofs of \eqref{1.2} and \eqref{1.3} are based upon an explicit enumeration of the elements of $\NR$.

Tripathi~\cite{Tripathi:1994} determined the set~$\cal{N}$ as
\begin{equation}
  {\cal N} = \left\{\vphantom\sum a\left(\Lfloor\frac{n-1}{s}\Rfloor+1\right) + dn,\ n=1,\ldots,a-1\right\},
  \label{eq:5.3.new}
\end{equation}
and used~\eqref{3.1} and~\eqref{3.2} to derive the above results of Roberts and Grant.

\subsection{The Sylvester sums}
It turns out that it is advantageous to use a slightly different, but equivalent formulation of~$\mathcal{N}$. 
Instead of~\eqref{5.3.new}, we use 
\begin{equation}
  {\cal N} = \left\{\vphantom\sum a\Lceil n/s\Rceil + dn,\ n=1,\ldots,a-1\right\}.
\end{equation}
This is easily seen as an equivalent by virtue of the 
relationship $\Lceil x/s\Rceil=\Lfloor(x-1)/s\Rfloor+1$, for all integers~$x$,
where $\Lceil x\Rceil$ denotes the least integer that is not less than~$x$.
Combining this with characterization~\eqref{3.1.3} gives:
\begin{equation}
  mS_{m-1}
  =
   a^{m-1}\sum_{n=0}^{a-1}B_m(\Lceil n/s\Rceil + nd/a)
  \,-\,B_m. 
 \label{eq:5.1.x.2}
\end{equation}
Note that for $s=1$, we recover formulae~\eqref{5.1.0} and~\eqref{4.1.x} from the two-variable case.
To arrive at a computational more convenient formulation, we determine the exponential generating function, 
and start with:
\begin{eqnarray*}
   \sum_{n\in{\cal N}_0} e^{nz}
  &=& 1+
   \sum_{n=1}^{a-1}e^{(a\Lceil n/s\Rceil+nd)z}  
   = 1+\sum_{n=1}^{a-1}e^{ndz} 
     \left[ 1+(e^{az}-1)\sum_{k=0}^{\Lceil n/s\Rceil-1} e^{kaz} \right] \\
  &=& 
    \frac{e^{adz}-1}{e^{dz}-1}
     + (e^{az}-1)   
      \sum_{n=1}^{a-1} \sum_{k=0}^{\Lceil n/s\Rceil-1} e^{(ka+nd)z}.  
\end{eqnarray*}
To evaluate the double sum, we change the order of summation: 
\begin{eqnarray*}
 \sum_{n=1}^{a-1} \sum_{k=0}^{\Lceil n/s\Rceil-1} e^{(ka+nd)z}
 &=& \sum_{k=0}^{\Lceil\frac{a-1}{s}\Rceil-1}     \sum_{n=ks+1}^{a-1} e^{(ka+nd)z} 
  = \sum_{k=0}^{\Lceil\frac{a-1}{s}\Rceil-1}e^{kaz}
     \left[\sum_{n=0}^{a-1} 
           e^{ndz}-
           \sum_{n=0}^{ks} e^{ndz}\right] \\ 
 &=& \sum_{k=0}^{\Lceil\frac{a-1}{s}\Rceil-1}e^{kaz}
     \left[\frac{e^{adz}         }{e^{dz}-1} -
           \frac{        e^{(ks+1)dz}}{e^{dz}-1}
     \right] \\ 
 &=& \frac{e^{\Lceil\frac{a-1}{s}\Rceil az}-1}{e^{az}-1}\frac{e^{adz}  }{e^{dz}-1}-
     \frac{e^{\Lceil\frac{a-1}{s}\Rceil bz}-1}{e^{bz}-1}\frac{e^{dz}}{e^{dz}-1}, 
\end{eqnarray*}
where we have set $b=a+sd$.
Combining all the elements, and simplifying the expression, gives
\begin{eqnarray}
 \sum_{n\in{\cal N}_0}e^{nz}
 &=&\frac{e^{\Lceil\frac{b-1}{s}\Rceil az}-1}
         {e^{dz}-1} +
    \left(e^{az}-1\right)\frac{e^{\Lceil\frac{a-1}{s}\Rceil
        bz}-1}{e^{bz}-1}\frac{1}{e^{-dz}-1}. 
 \label{eq:5.1.11}
\end{eqnarray}
Now use~\eqref{3.1.1} to give
\begin{eqnarray}
  \sum_{m=0}^\infty S_m\frac{z^m}{m!}
 &=&\frac{e^{\Lceil\frac{b-1}{s}\Rceil az}-1}{(e^{dz}-1)(e^{az}-1)} + 
    \frac{e^{\Lceil\frac{a-1}{s}\Rceil bz}-1}{(e^{-dz}-1)(e^{bz}-1)} - \frac{1}{e^z-1}. 
 \label{eq:5.4}
\end{eqnarray}
Note the simularity to~\eqref{5.x.1}, and, again, observe the footprint of the
Bernoulli polynomials.
As before, in order to invert this expression, we multiply both sides
of the equation by $z^2$,
and use the convolution property of the exponential generating
function to give
\begin{eqnarray}
  m(m-1)S_{m-2} 
 &=& \sum_{j=0}^m\BC{m}{j}d^{m-j-1}B_{m-j}a^{j-1}
     \phi_j\left(\Lceil(b-1)/s\Rceil\right) 
      + \nonumber \\
 & &
 \sum_{j=0}^m\BC{m}{j}(-d)^{m-j-1}B_{m-j}b^{j-1}
     \phi_j\left(\Lceil(a-1)/s\Rceil\right)
      - mB_{m-1}.
 \label{eq:5.7.y.1}
\end{eqnarray}
Taking account of the fact that the summands are zero for $j=0$, and expanding the Bernoullian polynomials, 
gives a generalisation of R\"odseth's formula:
\begin{eqnarray*}
  S_{m-1} 
 &=& \frac1{m(m+1)}\sum_{j=0}^{m}\sum_{i=0}^j\BC{m+1}{j+1}\BC{j+1}{i+1}d^{m-j-1}B_{m-j}B_{j-i}a^{j}\Lceil\frac{b-1}{s}\Rceil^{i+1} +  \\
 & & \frac1{m(m+1)}\sum_{j=0}^{m}\sum_{i=0}^j\BC{m+1}{j+1}\BC{j+1}{i+1}(-d)^{m-j-1}B_{m-j}B_{j-i}b^{j}\Lceil\frac{a-1}{s}\Rceil^{i+1} - \frac1mB_{m}.
\end{eqnarray*}
We note that this gives $S_m$ as a polynomial in $\Lceil(a-1)/s\Rceil$ of degree~$m+2$.
In particular, 
\begin{eqnarray*}
   2S_0 &=&
   -s\Lceil\frac{a-1}{s}\Rceil^2+(2a-2+s)\Lceil\frac{a-1}{s}\Rceil+(a-1)(d-1), \\
  12S_1 
  &=& -2s(2a+ds)\Lceil\frac{a-1}{s}\Rceil^3 + 
           \left(6a^2+3(2a+ds)(s-1)\right)\Lceil\frac{a-1}{s}\Rceil^2
           + \\
  & &  \left(6(d-1)a(a-1)+(3d-2a-ds)s\right)\Lceil\frac{a-1}{s}\Rceil   +
        (a-1)(d-1)(2ad-a-d-1), \\
  12S_2
  &=&  -s\left((a+sd)(2a+sd)+a^2\right)\Lceil\frac{a-1}{s}\Rceil^4 + \\
  & &  \left(4a^3+2(3s-2)a^2+2sd(3s-2)a+2s^2d^2(s-1)\right)\Lceil\frac{a-1}{s}\Rceil^3 + \\ 
  & &  \left(6(d-1)a^3+3(2-s-2d)a^2-3sd(s-2)a-sd^2(s^2-3s+1)\right)\Lceil\frac{a-1}{s}\Rceil^2 + \\
  & &  \left(2(2d-1)(d-1)a^3-(6d^2-6d+2)a^2+2d(d-s)a+s(1-s)d^2\right)\Lceil\frac{a-1}{s}\Rceil   + \\
  & &  ad(a-1)(d-1)(ad-a-d).
\end{eqnarray*}
The expressions for $S_1$ and $S_2$ in the arithmetic case have not previously appeared in the literature.
The first expression can be simplified to $S_0=(a-1+r)\Lceil(a-1)/s\Rceil + (a-1)(d-1)$, by repeatedly making the substitution $s\Lceil(a-1)/s\Rceil=a-1+s-r$, 
with $r=(a-1)\bmod s$, and is, of course, equivalent to~\eqref{1.3}.
Unfortunately, employing the same device for $S_1$ and $S_2$ does not
yield any expression that is much simpler than the ones given.
For the arithmetic case, the computationally more convenient form of
the Sylvester sums does not seem to have an 
apparent structure nor does it exhibit any pleasing pattern, as in the
two-variable case. At least, the author has not been able to find such. 
Nevertheless, one should bear in mind that~\eqref{5.7.y.1}, although not as aesthetically pleasing as~\eqref{5.1.x.2},
was derived from a purely computational point of view,
and for this purpose it is certainly adequate.

\section{The Frobenius problem for generalized arithmetic progressions}
Selmer~\cite{Selmer:1977} generalized the arithmetic progression by taking 
$\mathcal{A}=\left\{a, ha+d, ha+2d,\ldots,ha+sd\right\}$, with $d,h>0$ and
$\gcd(a,d)=1$.
As before, we may assume that $s<a$, as otherwise the last element can
be expressed as $ha+sd=d\times a + 1\times\left(ha+(s-a)d\right)$, and
thus the set $\NR$ does not alter if we remove the last element from~$\mathcal{A}$.
One can also show that, when $s<a$, there are no redundancies in 
$\mathcal{A}$, 
and none of its elements can be expressed as a nonnegative linear combination
of the others.
Selmer~\cite{Selmer:1977} showed that
\begin{equation}
  L=\left(h\Lfloor\frac{a-2}{s}\Rfloor+h-1\right)a+(a-1)d
\end{equation}
and
\begin{equation}
  S_0 = \frac12(a-1)\left(h\Lfloor\frac{a-1}{s}\Rfloor+d+h-1\right) + 
         \frac12h\left(\Lfloor\frac{a-1}{s}\Rfloor+1\right)(a-1\bmod s).
\end{equation}
Matthews~\cite{Matthews:2005} studied the generalized arithmetic sequence in the context of
numerical semigroups, and determined $\cal N$ as
\begin{equation}
 {\cal N} = \left\{ha\left(\Lfloor\frac{n-1}{s}\Rfloor+1\right)+dn,\ n=1,\ldots,a-1\right\}.
 \label{eq:6.0.3}
\end{equation}
She then applied~\eqref{3.1} and~\eqref{3.2} to derive Selmer's expressions for~$L$, and
the following expression for~$S_0$, that is easily seen to be equivalent to Selmer's:
\begin{equation}
   S_0 = 
   \frac12h\left(\Lfloor{a-2\over s}\Rfloor+1\right)(a+t)+\frac12(a-1)(d-1),
  \label{eq:6.1.3}
\end{equation}
where $t=(a-2)\bmod s$.

\subsection{The Sylvester sums}
As before, it is advantageous to rewrite~\eqref{6.0.3}, and take
\begin{equation}
 {\cal N} = \left\{ha\Lceil n/s \Rceil+dn,\ n=1,\ldots,a-1\right\}.
\end{equation} 
This allows us, as one would expect from the structure of $\mathcal{N}$,
to leverage the results from the previous section.
Characterization~\eqref{3.1.3} now gives:
\begin{equation}
  mS_{m-1}  =
   a^{m-1}\sum_{n=0}^{a-1}B_m(h\Lceil n/s\Rceil + nd/a)
  \,-\,B_m. 
 \label{eq:6.1.x.7}
\end{equation}
The computational more convenient form can be derived from the exponential generating function.
We start with
\begin{eqnarray*}
   \sum_{n\in{\cal N}_0} e^{nz}
  &=& \sum_{n=0}^{a-1}e^{(ha\Lceil n/s\Rceil+nd)z}
   =   \sum_{n=0}^{a-1}e^{(a\Lceil n/s\Rceil+nd/h)hz},
\end{eqnarray*}
and observe that this is the same sum as for the arithmetic case,
but with~$d/h$ and~$hz$, instead of~$d$ and~$z$. In the derivation of~\eqref{5.1.11}
we have not made use of the fact that~$d$ is an integer, so that we can use the result,
with the proviso that we replace~$z$ by~$hz$, and~$d$ by~$d/h$. This gives
\begin{eqnarray*}
   \sum_{n\in{\cal N}_0} e^{nz}
  &=& \frac{e^{\left(h\Lceil\frac{a-1}{s}\Rceil+d\right) az}-1}
         {e^{dz}-1} +
    \left(e^{ahz}-1\right)\frac{e^{\Lceil\frac{a-1}{s}\Rceil
        (ha+sd)z}-1}{e^{(ha+sd)z}-1}\frac{1}{e^{-dz}-1}. 
\end{eqnarray*}
Now use~\eqref{3.1.1} to give
\begin{eqnarray*}
  \sum_{m=0}^\infty S_m\frac{z^m}{m!}   
  &=& \frac{e^{\left(h\Lceil\frac{a-1}{s}\Rceil+d\right) az}-1}
         {(e^{dz}-1)(e^{az}-1)} +
    \frac{\left(e^{ahz}-1\right)\left(e^{\Lceil\frac{a-1}{s}\Rceil
        (ha+sd)z}-1\right)}{(e^{az}-1)(e^{(ha+sd)z}-1)(e^{-dz}-1)}\frac{}{} - \frac1{e^z-1}.
\end{eqnarray*}
To invert this expression, we multiply both sides by $z^3$, 
and use the convolution property of the exponential generating
function to give the companion of~\eqref{6.1.x.7}:
\begin{eqnarray}
 && \hspace{-1cm} m(m-1)(m-2)S_{m-3}  \nonumber \\ 
 &=& m\sum_{j=0}^{m-1}\BC{m-1}{j}d^{m-j-2}B_{m-j-1}a^{j-1}
     \phi_j\left(h\Lceil(a-1)/s\Rceil+d\right)      + \nonumber  \\
 & & \sum_{j=0}^m\sum_{i=0}^j\BC{m}{j}\BC{j}{i}a^{m-j-1}
     \phi_{m-j}(h)
     (-d)^{j-i-1}(ha+sd)^{i-1}B_{j-i}
      \phi_i\left(\Lceil(a-1)/s\Rceil\right) - \nonumber \\ 
 & & \vphantom{\sum_{j=0}^m} m(m-1)B_{m-2}.
\end{eqnarray}
Note that, for $h=1$, only the index $j=m-1$ of the first summation in the double sum remains, 
as $\phi_m(1)$ is zero for all~$m$, 
except for $m=1$, for which it has the value one. So that, for $h=1$, we recover the expression for the arithmetic case.
Expanding the Bernoullian polynomials does not give any additional
insights; we therefore refrain from doing this.

\section{Observations and comments}
We have taken $a$ as an arbitrary element of $\mathcal{A}$. 
It turns out that this choice, to a large extent, is irrelevant. 
Taking $f(n)=x^n$ in~\eqref{2.1} gives 
\[ \frac1{x^{a}-1}\sum_{n\in\mathcal{N}_0}x^n 
   = \frac1{x-1}+\sum_{n\in\NR}x^n,
\]
and shows that the left-hand side is an invariant with respect to the actual choice of~$a$. 
So, if we have $a_1$ and $a_2$, with corresponding sets $\mathcal{N}_0(a_1)$ and $\mathcal{N}_0(a_2)$, 
then knowledge of one set implies knowledge of the other through their generating functions:
\[ \frac1{x^{a_1}-1}\sum_{n\in\mathcal{N}_0(a_1)}x^n=\frac1{x^{a_2}-1}\sum_{n\in\mathcal{N}_0(a_2)}x^n.
\]
Therefore, the most reasonable choice seems to be to take $a$ as the smallest element of~$\mathcal{A}$,
as this minimizes the cardinality of the corresponding set $\mathcal{N}$. 
 
The special cases of the Frobenius problem that we considered here have elegant and succinct characterizations 
of the set~$\mathcal{N}$.
For the general case, however, it is not immediate what its structure is.
This is probably not an easy question to answer, given that just determining the largest element of ${\mathcal N}$
for variable~$k$ is NP-hard, as shown by Ram{\'\i}rez Alfons{\'\i}n~\cite{Ramirez:1996}.
However, for fixed~$k$, polynomial time algorithms do exist, see Kannan~\cite{Kannan:1992}, and algorithms to compute ${\mathcal N}$ have been given by 
Wilf~\cite{Wilf:1978}, and Nijenhuis~\cite{Nijenhuis:1979}, among others.

On a concluding, historical note, in~\cite{Tuenter:2006}, we showed that taking $f(n)=\Lfloor n/a\Rfloor$ in the two-variable case gives
\[ \sum_{n=1}^{a-1}\Lfloor nb/a\Rfloor = \sum_{n\in\NR}1=\frac12(a-1)(b-1),
\]
and provides a new interpretation of this sum as the number of integers that are not representable as a nonnegative linear
combination of~$a$ and~$b$, but omitted to mention that this formula can already be found in a paper by Hacks~\cite{Hacks:1893}.

\appendix

\section{Bernoulli numbers and polynomials}
In this section we have collected a number of properties of the 
Bernoulli numbers and polynomials in order not to distract from the flow
in the main body where these properties are needed.
The Bernoulli numbers can be defined by the implicit recurrence relation
\[ \sum_{j=0}^m\BC{m+1}{j}B_j=0, \ \ \ \mathrm{for} \ \ m\ge1,
\]
with $B_0=1$. This renders the first few Bernoulli numbers as $B_1=-1/2$, $B_2=1/6$, $B_3=0$, and $B_4=-1/30$.
Their exponential generating function is given by
\[ \sum_{m=0}^{\infty}B_m \frac{z^m}{m!}=\frac{z}{e^z-1}.
\]
The Bernoulli polynomials can be defined as
\[ B_m(x)=\sum_{r=0}^m\BC{m}{r}B_{m-r}x^r,
\]
to give $B_0(x)=1$, $B_1(x)=x-\frac12$, $B_2(x)=x^2-x+\frac16$,
$B_3(x)=x^3-\frac32x^2+\frac12x$, and $B_4(x)=x^4-2x^3+x^2-\frac1{30}$,
as the first few Bernoulli polynomials.
Their exponential generating function is given by 
\[ \sum_{m=0}^{\infty}B_m(x)\frac{z^m}{m!}=\frac{z}{e^z-1}e^{xz}.
\]
The Bernoulli polynomials satisfy the difference formula: $B_m(x+1)-B_m(x)=mx^{m-1}$,
and Raabe's multiplication theorem~\cite{Raabe:1851}: 
\[ \sum_{n=0}^{a-1}B_m(x+n/a)=a^{1-m}B_m(ax). \]
Both these are easily proved with the help of the exponential generating function.
A more detailed expos\'e on the Bernoulli numbers can be found in~\cite{GrahamKP:1989}.

\end{document}